%% file: PDF_Robustness_5.tex
\documentclass{siamonline1116}

\usepackage{braket,amsfonts}

\usepackage{array}

\usepackage[caption=false]{subfig}
\usepackage{pgfplots}

\newsiamthm{claim}{Claim}
\newsiamremark{rem}{Remark}
\newsiamremark{expl}{Example}
\newsiamremark{hypothesis}{Hypothesis}
\crefname{hypothesis}{Hypothesis}{Hypotheses}

\usepackage{algorithmic}

\usepackage{graphicx,epstopdf}

\Crefname{ALC@unique}{Line}{Lines}

\usepackage{amsopn}

\numberwithin{theorem}{section}

\usepackage{xspace}
\usepackage{bold-extra}
\usepackage[most]{tcolorbox}

\colorlet{texcscolor}{blue!50!black}
\colorlet{texemcolor}{red!70!black}
\colorlet{texpreamble}{red!70!black}
\colorlet{codebackground}{black!25!white!25}

\usepackage{amsmath,amssymb, amsfonts}
\usepackage{setspace}
\usepackage[mathscr]{euscript}
\usepackage{graphicx}
\usepackage{color}
\usepackage{cite}
\usepackage{geometry}
\usepackage{booktabs} 
\usepackage{enumerate}

\input{def}

\lstdefinestyle{siamlatex}{%
  style=tcblatex,
  texcsstyle=*\color{texcscolor},
  texcsstyle=[2]\color{texemcolor},
  keywordstyle=[2]\color{texemcolor},
  moretexcs={cref,Cref,maketitle,mathcal,text,headers,email,url},
}

\tcbset{%
  colframe=black!75!white!75,
  coltitle=white,
  colback=codebackground, 
  colbacklower=white, 
  fonttitle=\bfseries,
  arc=0pt,outer arc=0pt,
  top=1pt,bottom=1pt,left=1mm,right=1mm,middle=1mm,boxsep=1mm,
  leftrule=0.3mm,rightrule=0.3mm,toprule=0.3mm,bottomrule=0.3mm,
  listing options={style=siamlatex}
}

\newtcblisting[use counter=example]{example}[2][]{%
  title={Example~\thetcbcounter: #2},#1}

\newtcbinputlisting[use counter=example]{\examplefile}[3][]{%
  title={Example~\thetcbcounter: #2},listing file={#3},#1}

\DeclareTotalTCBox{\code}{ v O{} }
{ 
  fontupper=\ttfamily\color{black},
  nobeforeafter,
  tcbox raise base,
  colback=codebackground,colframe=white,
  top=0pt,bottom=0pt,left=0mm,right=0mm,
  leftrule=0pt,rightrule=0pt,toprule=0mm,bottomrule=0mm,
  boxsep=0.5mm,
  #2}{#1}

\patchcmd\newpage{\vfil}{}{}{}
\begin{tcbverbatimwrite}{tmp_\jobname_header.tex}
\title{Robustness of the Sobol' indices to distributional uncertainty
  \thanks{\funding{This work was supported by the National Science Foundation under grants DMS-1522765 and DMS-1745654.}}}

\author{Joseph Hart
  \thanks{Department of Mathematics, North Carolina State University, Raleigh, NC (\email{jlhart3@ncsu.edu}).}
  \and
  Pierre Gremaud
   \thanks{Department of Mathematics, North Carolina State University, Raleigh, NC (\email{gremaud@ncsu.edu}).}
}

\headers{Robustness of Sobol' indices}
{Joseph Hart and Pierre Gremaud}
\end{tcbverbatimwrite}
\input{tmp_\jobname_header.tex}

\begin{document}
\maketitle

\begin{tcbverbatimwrite}{tmp_\jobname_abstract.tex}
\begin{abstract}
Global sensitivity analysis (GSA) is used to quantify the influence of uncertain variables in a mathematical model. Prior to performing GSA, the user must specify (or implicitly assume), a probability distribution to model the uncertainty, and possibly statistical dependencies, of the variables. Determining this distribution is challenging in practice as the user has limited and imprecise knowledge of the uncertain variables. This article analyzes the robustness of the Sobol' indices, a commonly used tool in GSA, to changes in the distribution of the uncertain variables. A method for assessing such robustness is developed which requires minimal user specification and no additional evaluations of the model. Theoretical and computational aspects of the method are considered and illustrated through examples.
\end{abstract}


\end{tcbverbatimwrite}
\input{tmp_\jobname_abstract.tex}

\section{Introduction}
Global sensitivity analysis (GSA) aims to quantify the relative importance of the input variables $\X=(X_1,X_2, \dots, X_p)$ to a function $f$ \cite{saltellibook}. Such quantification is a crucial step in the development of predictive models. The Sobol' indices \cite{sobol93,sobol} are a commonly used tool for this task, though there are many alternative methods as well, see \cite{saltellibook,iooss,uq_handbook} and references therein for overviews of GSA. Performing GSA in practice consists of:
\begin{enumerate}
\item[i.] defining a probability distribution for $\X$,
\item[ii.] evaluating $f$ at samples from this distribution,
\item[iii.] and computing statistics using these evaluations.
\end{enumerate}
 Much of the GSA research has focused on steps (ii) and (iii) assuming that step (i) is done by the user. In practice, defining a probability distribution for $\X$ is a challenging task which carries uncertainties itself. This raises the question, ``how robust is my GSA to changes in the distribution of $\X$?" In this article we address this question for the Sobol' indices. 

The robustness of the Sobol' indices with respect to changes in the marginal distributions of $X_k$, $k=1,2,\dots,p$, is considered in the Life Cycle Analysis literature \cite{lca_app}. Their approach requires the user to specify various admissible marginal distributions for the variables and compute the Sobol' indices for each possibility. This requires many evaluations of $f$ thus limiting its use more broadly. The robustness of the Sobol' indices with respect to changes in the correlation structure is highlighted in \cite{lca_correlations} where the authors seek to quantify the risk of ignoring correlations between the variables. Imprecise probabilities are used in \cite{hall} to quantify uncertainty in the sensitivity indices. This approach requires the user to parameterize admissible distributions and collect additional samples, i.e. additional evaluations of $f$. In the ecology literature \cite{WARM}, the robustness of Sobol' indices to changes in the means and standard deviations defining normally distributed inputs is examined. The challenges of deep uncertainty, i.e. uncertainty in the distributions of $X_k$, $k=1,2,\dots,p$, are identified in \cite{deep_uncertainty} where the authors compute Sobol' indices for different distributions and define ``robust sensitivity indicators" as a function of the Sobol' indices from different distributions. Changes in the marginal distributions were shown to change the ordering of the Sobol' indices in \cite{cousins}. Similar questions about the robustness of computed quantities with respect to distributional uncertainty may be found in \cite{climate_app,dice,chick,beckman_mckay}.

In \cite{anova_mult_dist}, the functional analysis of variance decomposition (the building block of Sobol' indices) is analyzed when $\X$ does not have a unique distribution but rather multiple possible distributions. The authors provide a framework for analyzing the robustness of the Sobol' indices which depends upon the user specifying a prior on the space of possible distributions of $\X$. In line with \cite{SA_review_techreport}, the robustness of the Sobol' indices may be determined by considering a set of possible distributions, sampling from their mixture distribution, and computing the Sobol' indices with respect to each distribution using a weighting scheme. This approach does not require additional evaluations of $f$, but the user must specify the set of possible distributions, which is challenging in practice.

All of the aforementioned approaches require user specification of possible distributions, additional evaluations of $f$, or both. In this article, we present a method to measure the robustness of the Sobol' indices to distributional uncertainties without requiring either. In particular, we consider perturbations of the probability density function (PDF) of $\X$ and compute the extreme scenarios when the Sobol' indices differ most from those computed with the user specified PDF. A judicious formulation allows us to determine these extreme scenarios as the solution of an optimization problem which is solved in closed form. The Sobol' indices with a perturbed PDF are then computed using weighted averages. Our proposed method is a post processing step which requires minimal user specification and no evaluations of $f$ beyond those already taken to compute the Sobol' indices.

Section~\ref{sec:sobol_review} provides a review of the Sobol' indices and establishes notation for the article. Section~\ref{sec:general} develops the theoretical and computational framework for our method of assessing robustness. In some cases, the Sobol' indices may not be robust; however, the importance of the input variables relative to one another may be. This motives us to define and study ``normalized Sobol' indices" in Section~\ref{sec:normalized}. In Section~\ref{sec:algorithm}, we provide an algorithmic description of our method and give a guide for visualizing and interpreting the results. A variety of examples are given in Section~\ref{sec:numerical_results} to demonstrate the proposed approach and highlight its properties. Section~\ref{sec:conclusion} concludes the article by elaborating on limitations and possible extensions of this work.

\section{Review of Sobol' Indices}
\label{sec:sobol_review}

Let $f:\Omega \to \R$, $\Omega = \Omega_1 \times \Omega_2 \times \cdot \cdot \cdot \times \Omega_p \subset \R^p$, be a function or model and  let $\X=(X_1,X_2,\dots,X_p) \in \Omega$ be the input variables of that model. The Sobol' indices \cite{sobol93,sobol}  measure the importance of a variable (or group of variables) by apportioning to the variable (or group of variables) its relative contribution to the variance of $f(\X)$. 

Let $u=\{i_1,i_2,\dots,i_k\}$ be a  subset of $\{1,2,\dots ,p\}$ and $\sim u=\{1,2,\dots,p\}\setminus u$ be its complement. We refer to the group of variables corresponding to $u$ as  $\X_u=(X_{i_1},X_{i_2},\dots,X_{i_k})$. 
Following \cite{kucherenko}, assume that $f(\X)$ is square integrable and consider the law of total variance decomposition,
\begin{eqnarray}
\label{law_total_variance}
\Va(f(\X)) = \Va(\E[f(\X) \vert \X_u]) + \E[\Va(f(\X) \vert \X_u)].
\end{eqnarray}
Using \eqref{law_total_variance}, the Sobol' index and total Sobol' index for $\X_u$ are defined as
\begin{eqnarray*}
S_u = \frac{\Va(\E[f(\X) \vert \X_u]) }{\Va(f(\X))} \qquad \text{and} \qquad T_u = 1-S_{\sim u},
\end{eqnarray*}
respectively. The Sobol' index $S_u$ may be interpreted as the proportion of $\Va(f(\X))$ contributed by $\X_u$ alone; hence $S_u \in [0,1]$ and larger values indicate that $\X_u$ is influential. The total Sobol' index $T_u$ may be interpreted as the proportion of $\Va(f(\X))$ remaining if $\X_{\sim u}$ is known, hence $T_u \in [0,1]$ and larger values indicate that $\X_u$ is influential. Under the assumption that $X_1,X_2,\dots,X_p$ are independent, $S_u \le T_u$ and their difference may be interpreted as a measure of the interaction between $\X_u$ and $\X_{\sim u}$. They possess other useful statistical properties and are a preferred method for global sensitivity analysis in many applications. However, much of the statistical theory does not generalize when $X_1,X_2,\dots,X_p$ possess dependencies. In \cite{hart_corr_var}, $T_u$ is shown to have a useful approximation theoretic interpretation (with both independent or dependent inputs). Specifically, $T_u$ corresponds to the squared relative $L^2(\Omega)$ error when $f(\X)$ is optimally approximated by a function which does not depend on $\X_u$. In other words, the influence of $\X_u$ is measured by the error when $f(\X)$ is approximated by a function which does not depend on $\X_u$. Using this interpretation, the total Sobol' index provides a useful measure of the importance of $\X_u$, with independent or dependent variables. Statistical dependencies in $\X$ may effect the magnitude of $T_u$ and care must be taken in how one measures the relative importance of the variables. In Section~\ref{sec:normalized}, we introduce the \textit{normalized Sobol' index} as a means to measure the relative importance of variables when dependencies exist.

We direct the reader to \cite{sobol93,sobol,iooss,hart_corr_var,Sobol_UQ_handbook,sobol2003,saltelli2010,kucherenko,fast_dependent_variables,mara} for a fuller discussion of the Sobol' indices, total Sobol' indices, their interpretations, and their estimation. The reader may also consider \cite{iooss_prieur_shapley,owen_prieur_shapley,staum} for additional discussion of Sobol' indices with dependent variables, and an alternative approach, the Shapley value.

\section{Robustness of the Sobol' Index to PDF Perturbations}
\label{sec:general}
Assume that $\X$ admits a PDF $\phi$. For simplicity, and because of its approximation theoretic interpretation with dependent variables, we focus on the robustness of the total Sobol' index $T_u$ to changes in $\phi$; the Sobol' index $S_u$ may be analyzed in a similar fashion. For the remainder of the article, we will use ``Sobol' index" to refer to the total Sobol' index. 

There are multiple ways to express and estimate $T_u$; a useful expression from \cite{kucherenko} is
\begin{eqnarray}
\label{conditional_estimator}
T_u = \frac{ \frac{1}{2} \int_{\Omega \times \Omega_u} (f(\x) - f(\x') )^2  \phi(\x) \phi_{\x \vert \x_{\sim u}} (\x' \vert \x_{\sim u}) \,d\x \,d\x_u' }{\int_{\Omega} f(\x)^2 \phi(\x) d\x - \left( \int_{\Omega} f(\x) \phi(\x) d\x \right)^2}
\end{eqnarray}
where $\x=(\x_u,\x_{\sim u}$), $\x'=(\x_u',\x_{\sim u}$), $\phi_{\x \vert \x_{\sim u}}$ is the conditional density for $\X \vert \X_{\sim u}$, and $\Omega_u$ is the Cartesian product of each $\Omega_k$, $k \in u$. Note that $\x=(\x_u,\x_{\sim u}$) is not a permutation of the entries of $\x$ but rather a partitioning of them. Then $T_u$ may be estimated by drawing samples from $\X$, drawing a second set of samples from $\X \vert \X_{\sim u}$, and estimating \eqref{conditional_estimator} via Monte Carlo integration of the numerator and denominator separately.

The basic idea of the proposed approach is to view $T_u$ as an operator which inputs the PDF and returns the Sobol' index. We compute the Fr\'echet derivative of this operator at $\phi$ and use it to analyze the robustness of $T_u$. To this end, we make the following assumptions throughout the article:
\begin{enumerate}
\item $\Omega$ is a Cartesian product of compact intervals,
\item $\phi(\x) > 0$ $\forall \x \in \Omega$,
\item $\phi$ is continuous on $\Omega$,
\item $f$ is bounded on $\Omega$.
\end{enumerate}
Some of the results below may be shown with weaker assumptions, these overarching assumptions are made now for conciseness and simplicity. Without loss of generality, under the assumptions above, assume that $\Omega = [0,1]^p$. We revisit our assumption on the compactness of $\Omega$ in Section~\ref{sec:conclusion}.

We seek to perturb the PDF, so it is essential that the perturbations preserve properties of PDF's, specifically, that every PDF is non negative and its integral over $\Omega$ equals one.

Since $\phi>0$ is continuous and $\Omega$ is compact, $\phi$ is bounded above and below by positive real numbers. Define the Banach space $V$ as the set of all bounded functions on $\Omega$ equipped with the norm
\begin{eqnarray*}
\vert \vert \psi \vert \vert_V = \left\vert \left\vert \frac{\psi}{\phi} \right\vert \right\vert_{L^\infty(\Omega)},
\end{eqnarray*}
where $\vert \vert \cdot \vert \vert_{L^\infty(\Omega)}$ is the supremum norm on $L^\infty(\Omega)$, the set of bounded function on $\Omega$. This norm ensures that $\phi+\psi \ge 0$ for every $\psi \in V$ with $\vert \vert \psi \vert \vert_V \le 1$, the non negativity property of PDF's.

To ensure that the integral over $\Omega$ equals one, we introduce a normalization operator which takes $\eta \in V$ and returns $\frac{\eta}{\int_\Omega \eta(\x)d\x}$. Composing this normalization operator with \eqref{conditional_estimator} yields the Sobol' index as an operator on $V$. Define $F,G,T_u:V \to \R$ by
\begin{eqnarray}
\label{F}
F(\eta) = \frac{1}{2} \int_{\Omega \times \Omega_u} (f(\x)-f(\x'))^2 \eta(\x) \eta(\x') \frac{1}{\int_{\Omega_u} \eta(\x)d\x_u} d\x d\x_u',
\end{eqnarray}
 \begin{eqnarray}
 \label{G}
G(\eta) = \int_{\Omega} f(\x)^2 \eta(\x) d\x - \frac{1}{\int_\Omega \eta(\x)d\x} \left( \int_{\Omega} f(\x) \eta(\x) d\x \right)^2,
\end{eqnarray}
\begin{eqnarray}
\label{T_u}
T_u(\eta) = \frac{F(\eta)}{G(\eta)} .
\end{eqnarray}
It is easily observed that multiplying the numerator and denominator of  \eqref{T_u} by $ \frac{1}{\int_\Omega \eta(\x)d\x}$ yields that \eqref{T_u} and  \eqref{conditional_estimator} coincide with $\phi$ replaced by $ \frac{\eta}{\int_\Omega \eta(\x)d\x}$. In this framework, every $\eta \in V$ such that $\vert \vert \phi-\eta \vert \vert_V \le 1$ is nonnegative and $T_u(\eta)$ corresponds to the Sobol' index computed with respect to the PDF $\frac{\eta}{\int_\Omega \eta(\x)d\x}$.

Having defined the Sobol' index as an operator which inputs bounded PDF's, Theorem~\ref{thm:frechet} below gives the Fr\'echet derivative of the Sobol' index at $\phi$.

\begin{theorem}
\label{thm:frechet}
The operator $T_u$ is Fr\'echet differentiable at $\phi$ with Fr\'echet derivative $\D T_u(\phi): V \to \R$ given by the bounded linear operator
\begin{eqnarray}
\D T_u(\phi) \psi = \frac{\D F(\phi)\psi}{G(\phi)} - T_u(\phi) \frac{\D G(\phi)\psi}{G(\phi)},
\end{eqnarray}
where
\begin{align*}
\D F(\phi)\psi =& \frac{1}{2} \int_{\Omega \times \Omega_u} (f(\x)-f(\x'))^2  \frac{\psi(\x')}{\phi(\x')} \phi(\x) \phi_{\x \vert \x_{\sim u}} (\x' \vert \x_{\sim u})  d\x d\x_u' \\
+& \frac{1}{2} \int_{\Omega \times \Omega_u} (f(\x)-f(\x'))^2 \frac{\psi(\x)}{\phi(\x)} \phi(\x) \phi_{\x \vert \x_{\sim u}} (\x' \vert \x_{\sim u})  d\x d\x_u' \\
-& \frac{1}{2} \int_{\Omega \times \Omega_u} (f(\x)-f(\x'))^2 \frac{\int_{\Omega_u} \psi(\x) d\x_u}{\int_{\Omega_u} \phi(\x) d\x_u} \phi(\x) \phi_{\x \vert \x_{\sim u}} (\x' \vert \x_{\sim u}) d\x d\x_u' \\
\end{align*}
and
\begin{align*}
\D G(\phi)\psi = & \int_\Omega f(\x)^2 \frac{\psi(\x)}{\phi(\x)} \phi(\x) d\x \\
& - 2 \int_\Omega f(\x)\phi(\x) d\x \int_\Omega f(\x) \frac{\psi(\x)}{\phi(\x)} \phi(\x)d\x \\
& + \left(\int_\Omega \frac{\psi(\x)}{\phi(\x)} \phi(\x)d\x \right) \left( \int_\Omega f(\x) \phi(\x)d\x \right)^2 .
\end{align*}
\end{theorem}

A proof for Theorem~\ref{thm:frechet} is given in the appendix. Without the assumption that $\Omega$ is compact (or at least bounded), an infinitesimal perturbation $\epsilon \psi$ may yield that $\phi+\epsilon \psi$ is not integrable. Hence it is a theoretical necessity to assume that $\Omega$ is bounded. In Section~\ref{sec:conclusion}, we revisit this discussion and highlight that this assumption is less restrictive in practice.

If the Sobol' index is computed using Monte Carlo estimators of \eqref{conditional_estimator}, then $\D T_u(\phi)\psi$ may be estimated using these samples and evaluations of $f$; the only additional work is evaluating $\phi$ and $\psi$ at the sample points. Hence $\D T_u(\phi)\psi$ may be estimated at any $\psi \in V$ with negligible computational cost. This is why, as previously mentioned, our method is a post processing step which requires no additional evaluations of $f$ beyond those taken to compute the Sobol' indices.

We seek to find an ``optimal" perturbation of $\phi$ in the sense that it causes the greatest change in the Sobol' index. The locally optimal perturbation is the $\psi \in V$, $\vert \vert \psi \vert \vert_V \le 1$, which maximizes $\vert \D T_u(\phi)\psi \vert$. To estimate this $\psi$, we define a finite dimensional subspace $V_M \subset V$ and compute the operator norm of the restriction of $\D T_u(\phi)$ to $V_M$. When choosing $V_M$, there is a trade off to consider between the approximating properties of functions from $V_M$, our ability to use existing samples to estimate the action of $\D T_u(\phi)$ on functions from $V_M$, and the ease of computing the operator norm of $\D T_u(\phi)$ restricted to $V_M$. In what follows, we choose $V_M$ to be a subspace generated by the span of a set of locally supported piecewise constant functions. 

Let $R_i$, $i=1,2,\dots,M$, be a partition of $\Omega$ into open hyperrectangles, i.e. $\Omega = \cup_{i=1}^M \overline{R_i}$ and $R_i \cap R_j = \emptyset$ for $i \ne j$; $\overline{R_i}$ denotes the closure of $R_i$. Define 
\[ \psi_i(\x) = \begin{cases} 
      1 & \x \in R_i \\
      0 & \x \notin R_i\\
   \end{cases}
\]
to be the indicator function of $R_i$, $i=1,2,\dots,M$, and $V_M=span\{\psi_1,\psi_2,\dots,\psi_M\}$, a $M$ dimensional subspace of $V$. The partition may be efficiently constructed using Regression Trees \cite{breimanbook}; we will elaborate on this in Section~\ref{sec:algorithm}. 

Constructing $V_M$ in this way has computational and approximation theoretic advantages. Its computational advantage, demonstrated below, is that it enables a closed-form solution to an otherwise challenging optimization problem. Its approximation theoretic advantage is that is provides a mechanism to constrain the subspace $V_M$ with the existing evaluations of $f$, see Section~\ref{sec:algorithm} for more details. Piecewise constant functions defined on a partition of the domain are useful in the sense that they permit perturbations of any functional form, constrained by the coarseness of the partition. The coarseness of the partition depends on the number of existing evaluations of $f$, hence the form of functions from $V_M$ are constrained by the existing data, not the users imposition of functional forms. The proposed subspace $V_M$ is optimal in the (informal) sense that is provides the most flexibility in functional forms given the constraint of existing data.

The operator norm of $\D T_u(\phi)$ restricted to $V_M$ is given by
\begin{align*}
\vert \vert \D T_u(\phi) \vert \vert_{\LL(V_M,\R)} & = \max_{\substack{\psi \in V_M \\ \vert \vert \psi \vert \vert_V \le 1 }} \vert \D T_u(\phi)\psi \vert \\
& = \max_{\substack{\mathbf a \in \R^M\\ \vert \vert \sum_{i=1}^M a_i \psi_i \vert \vert_V \le 1}} \left\vert \D T_u(\phi) \left(\sum_{i=1}^M a_i \psi_i \right) \right\vert \\
& = \max_{\substack{\mathbf a \in \R^M\\ \vert \vert \sum_{i=1}^M a_i \psi_i \vert \vert_V \le 1}} \left\vert \sum_{i=1}^M a_i \D T_u(\phi) \psi_i \right\vert 
\end{align*}

Since the basis functions have disjoint support, it follows that
\begin{eqnarray*}
\left\vert \left\vert \sum_{i=1}^M a_i \psi_i \right\vert \right\vert_V =  \left\vert \left\vert \sum_{i=1}^M a_i \frac{1}{\phi} \psi_i  \right\vert \right\vert_{L^\infty(\Omega)} = \max\limits_{i=1,2,\dots,M} \vert a_i \vert \left\vert \left\vert \frac{1}{\phi}  \right\vert \right\vert_{L^\infty(R_i)} ,
\end{eqnarray*}
which implies
\begin{eqnarray*}
\left\vert \left\vert \sum_{i=1}^M a_i \psi_i \right\vert \right\vert_V  \le 1
\end{eqnarray*}
is equivalent to $\vert a_i \vert \le b_i$, where $b_i$ is the infimum of $\phi$ on $R_i$, $i=1,2,\dots,M$.

Let $\mathbf d \in \R^M$ be defined by $d_i = \D T_u(\phi)\psi_i$, $i=1,2,\dots,M$. Then we have
\begin{eqnarray*}
\vert \vert \D T_u(\phi) \vert \vert_{\LL(V_M,\R)} =  \max_{\substack{\mathbf a \in \R^M \\  \vert a_i \vert \le b_i \\ i=1,2,\dots,M }} \vert \mathbf d^T \mathbf a \vert  .
\end{eqnarray*}

This problem may be solved in closed form to get 
\begin{eqnarray*}
a_i = \text{sign}(d_i) b_i
\end{eqnarray*}
 and
\begin{eqnarray*}
\vert \vert \D T_u(\phi) \vert \vert_{\LL(V_M,\R)} = \vert \vert \mathbf d \vert \vert_1 .
\end{eqnarray*}

In what follows, we refer to $\psi \in V_M$, $\vert \vert \psi \vert \vert_V \le 1$, which maximizes the Fr\'echet derivative, as the \textit{optimal perturbation}. Finding the optimal perturbation and the corresponding operator norm simplifies to evaluating $\D T_u(\phi)\psi_i$ for $i=1,2,\dots,M$, which may be estimated with negligible additional computation. However, estimating $\D T_u(\phi)\psi_i$ is typically more challenging than estimating the Sobol' index. Rather than inferring robustness with $\vert \vert \D T_u(\phi) \vert \vert_{\LL(V_M,\R)}$, we propose to:
\begin{enumerate}
\item[i.] estimate $a_i = \text{sign}(d_i) b_i$, $i=1,2,\dots,M$,
\item[ii.] use weighted averaging with the existing evaluations of $f$ and $\phi$ to estimate the Sobol' indices with respect to the \textit{optimally perturbed PDF}, which we define as
\end{enumerate}
\begin{eqnarray}
\label{perturbed_pdf}
\frac{\phi+\delta \sum\limits_{i=1}^M a_i \psi_i}{1+\delta \sum\limits_{i=1}^M a_i \text{vol}(R_i)},
\end{eqnarray}
where $\delta \in [-1,1]$ is a parameter to scale the size of the perturbation and $\text{vol}(R_i)$ is the volume of the set $R_i$; the determination of $\delta$ will be discussed in Section~\ref{sec:algorithm}. We will refer to the Sobol' indices computed with $\phi$ as the \textit{nominal Sobol' indices} and the Sobol' indices computed with the optimally perturbed PDF \eqref{perturbed_pdf} as the \textit{perturbed Sobol' indices}.

In practice, it is suggested to estimate the terms $\text{vol}(R_i)$, $i=1,2,\dots,M$, in \eqref{perturbed_pdf} with a Monte Carlo estimator from the existing data. They may be computed analytically since $R_i$ is known; however, if they are computed exactly then the weights used to estimate perturbed Sobol' indices may not sum to one because of Monte Carlo error in the estimate. This can bias the resulting analysis. Estimating $\text{vol}(R_i)$, $i=1,2,\dots,M$, from the existing data diminishes this potential bias.

Our weighted averaging approach is an improvement from traditional derivative based robustness analysis in several ways:
\begin{enumerate}
\item[$\bullet$] Estimating $a_i$ is easy. Since $\phi$ is known, $b_i$ may be computed numerically by querying the existing evaluations of $\phi$ (or possibly analytically, for instance if $\phi \equiv 1$ then $b_i=1$ for every $i$), we consider this negligible. Assuming that we have enough samples for the Sobol' index estimation to converge, determining the sign of $d_i$ with these samples is relatively easy. Additionally, when we do not determine the sign of $d_i$ correctly it is frequently because $\D T_u(\phi)\psi_i \approx 0$, in which case this error is benign in the scope of our analysis.
\item[$\bullet$] The user must determine $\delta$; however, various values of $\delta \in [-1,1]$ may be tested at negligible computation cost. The sample standard deviation of the weighted average may be compared with the sample standard deviation in the original estimator to determine admissible values of $\delta$. Additional details are given in Section~\ref{sec:algorithm}.
\item[$\bullet$] Computing the perturbed Sobol' indices estimates a realized worst case. This is superior to worst case bounds, error bars, or confidence intervals, which in many cases are overly pessimistic. Further, computing error bars for each Sobol' index individually may yield misleading results. For instance, error bars for two variables may yield large intervals for each Sobol' index, but their magnitude relative to one another is nearly constant for any PDF perturbation. In this case the user would incorrectly conclude that the relative importance of the variables to one another is uncertain.

\end{enumerate}

As previously highlighted, one way to test for robustness is to use weighted averages to estimate the Sobol' indices with different PDF's. The challenge with this approach is that the user must specify the perturbed PDF's. Our method may be viewed as an improvement on this idea by automating the choice of perturbed PDF's. The Fr\'echet derivative operator norm yields a locally optimal perturbation, which will likely reveal greater changes in the Sobol' indices when compared with a user manually selecting a small set of perturbed PDF's. However, our method does not have the danger of finding unrealistic worst cases since it only seeks perturbations in a neighborhood of the existing PDF and is constrained to use the existing samples.

\section{Robustness of the Normalized Sobol' Index to PDF Perturbations}
\label{sec:normalized}

As highlighted in \cite{hart_corr_var}, the Sobol' indices are frequently smaller when $\X$ possesses stronger dependencies. Since the magnitude of the Sobol' indices may change when the PDF is perturbed, we seek to analyze the robustness of the relative importance of the variables rather than the magnitude of the Sobol' indices. This is useful in practice when, for instance, a modeler (for lack of better knowledge) assumes the variables are independent, computes the Sobol' indices to rank the importance of the variables, but is uncertain of their ranking because dependencies that were potentially ignored. For clarity and notational simplicity, the remainder of the article will focus on the Sobol' indices when $u=\{k\}$ is a singleton, i.e. the set of Sobol' indices $\{T_k\}_{k=1}^p$. 

To measure the relative importance of the variables as the PDF varies, define the normalized Sobol' index $\T_k:V \to \R$ as
\begin{eqnarray}
\label{normized_sobol}
\T_k(\phi) = \frac{T_k(\phi)}{\sum\limits_{i=1}^p T_i(\phi)} 
\end{eqnarray}
for $k=1,2,\dots,p$. The example below illustrates the behavior of the Sobol' indices and normalized Sobol' indices.

\subsection*{Example}
\begin{em}
Let 
\begin{eqnarray}
f(\X)=1.5X_1+1.25X_2+X_3\label{simpleex}
\end{eqnarray}
and $\X$ follow a multivariate normal distribution with mean $\mu$ and covariance matrix $\Sigma$ given by
\[ \mathbf{\mu}= \left[ \begin{array}{cc}
0 \\
0 \\
0 \\
\end{array} \right],\hspace{10 mm}
\Sigma= \left[ \begin{array}{ccccc}
1 & \rho & \rho \\
\rho & 1 & \rho \\
\rho & \rho & 1\\
\end{array} \right], \qquad 0\le \rho\le 1. \] 

Figure~\ref{fig:corrknob} shows the Sobol' indices $T_k$ and normalized Sobol' indices $\T_k$, $k=1,2,3$, as a functions of $\rho$. As $\rho$ increases (the correlations strengthen), the Sobol' indices decrease while the normalized Sobol' indices are constant. The trend of $T_k$ decreasing as correlations increase is general; \cite{hart_corr_var} explains it with an approximation theoretic perspective of the Sobol' indices. The normalized Sobol' indices are constant indicating that, though the Sobol' indices decrease, the relative importance of the variables does not change. In general the normalized Sobol' indices will change as the distribution of the input variables changes.

\begin{figure}[h]
\centering
\includegraphics[width=.45 \textwidth]{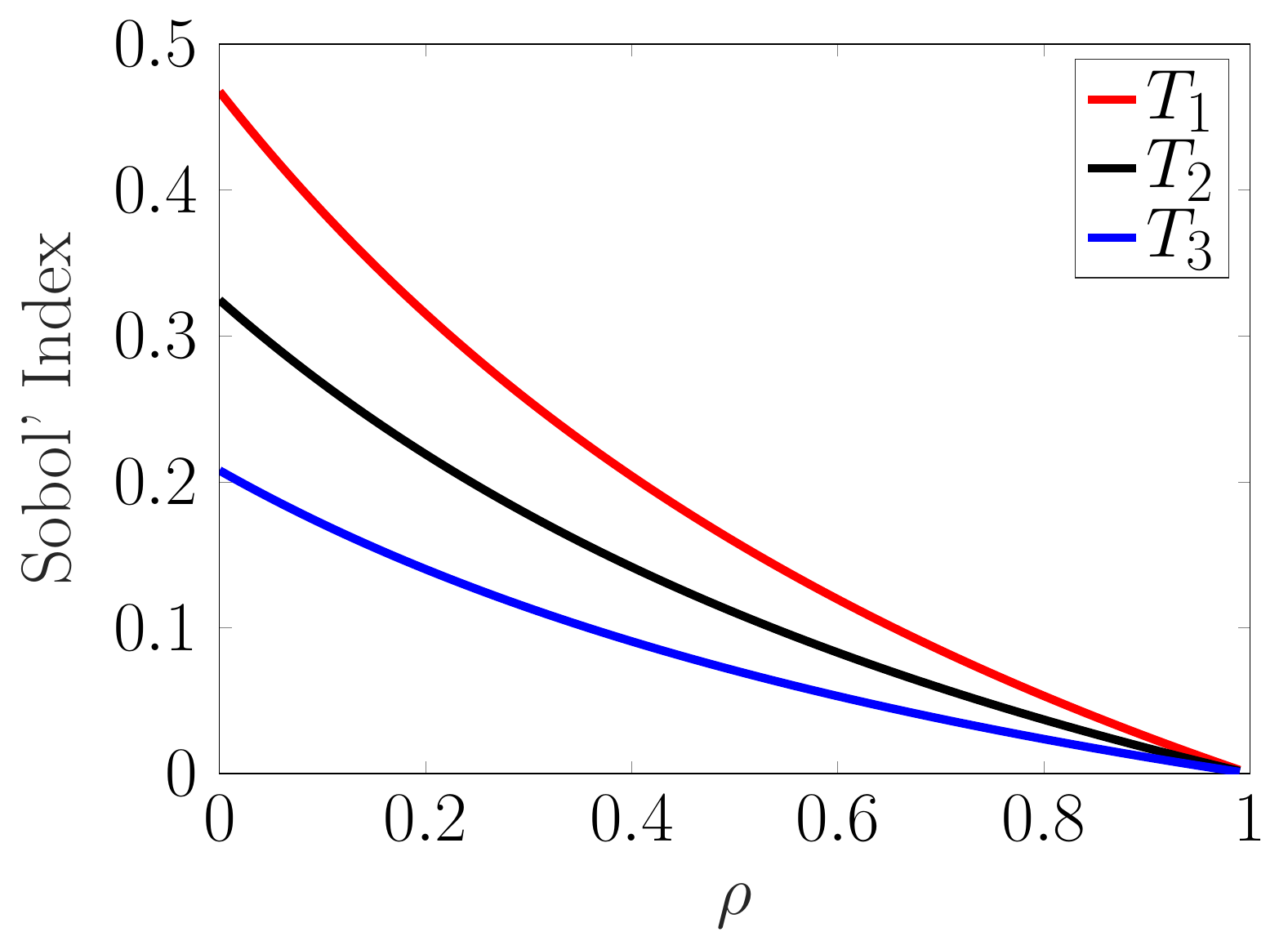}
\includegraphics[width=.45 \textwidth]{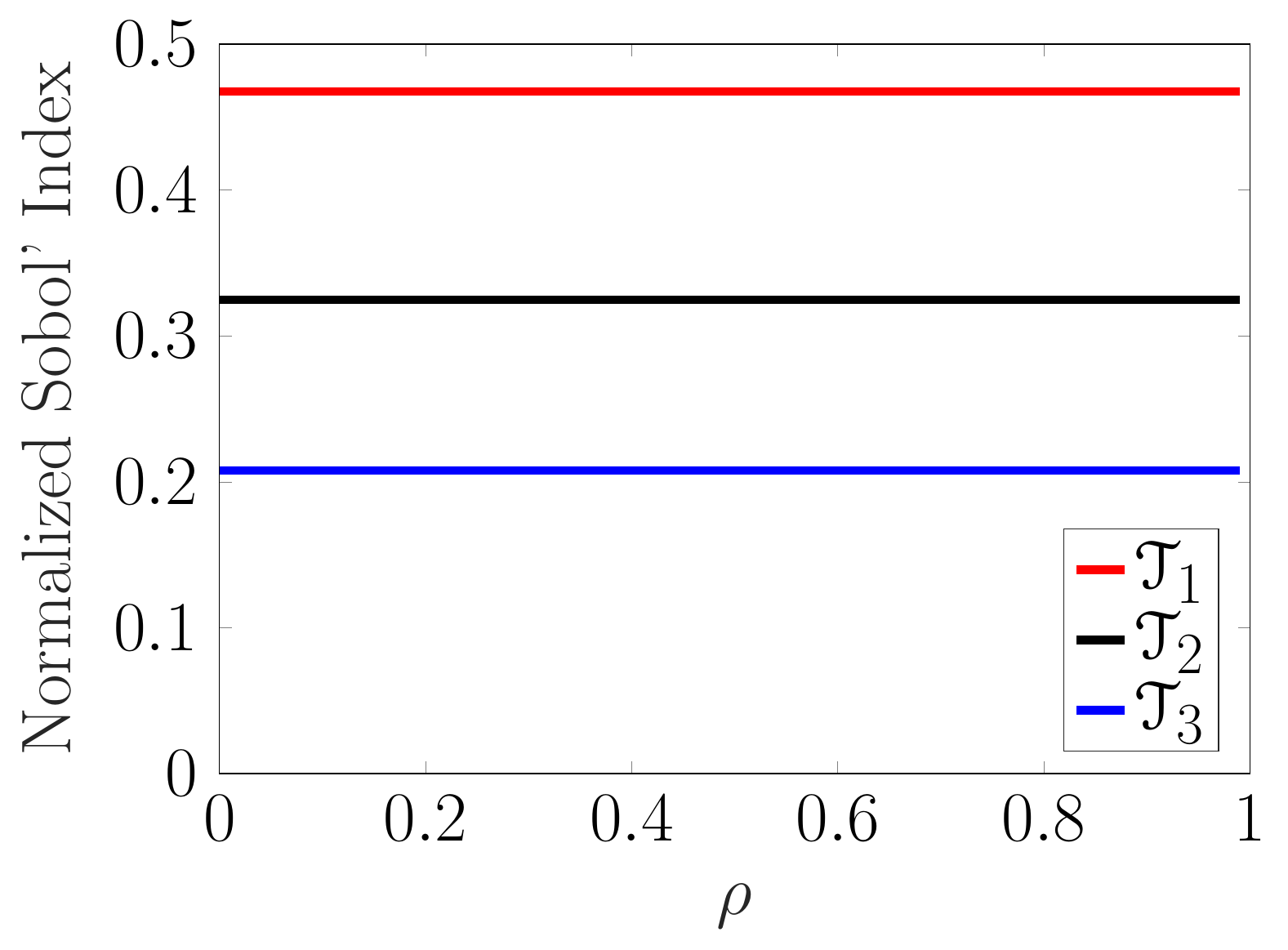}
\caption{Sobol' indices (left) and normalized Sobol' indices (right) for \eqref{simpleex} with increasing correlation strength as $\rho$ varies from 0 to 1.}
\label{fig:corrknob}
\end{figure}

\end{em}

Applying Theorem~\ref{thm:frechet} and the quotient rule to \eqref{normized_sobol} yields that $\T_k$ is Fr\'echet differentiable with Fr\'echet derivative
\begin{eqnarray*}
\D \T_k(\phi)\psi = \frac{ \left( \sum\limits_{i=1}^p T_i(\phi) \right) \D T_k(\phi)\psi - T_k(\phi) \left( \sum\limits_{i=1}^p \D T_i(\phi)\psi \right)}{ \left( \sum\limits_{i=1}^p T_i(\phi) \right)^2} .
\end{eqnarray*}

Since $\D \T_k(\phi)$ is a linear combination of the operators $\D T_k(\phi)$ from Section~\ref{sec:general}, we may easily estimate $\D \T_k(\phi) \psi$ using the same results previously presented. In fact, in Section~\ref{sec:general} a subspace $V_M = span\{\psi_1,\psi_2,\dots,\psi_M\}$ is defined and we compute $\D T_k(\phi)\psi_i$ for $i=1,2,\dots,M$. Using this computation we may easily compute $\D \T_k(\phi)\psi_i$, $i=1,2,\dots,M$, at no additional cost. The same procedure from Section~\ref{sec:general} may be adopted to compute perturbed PDF's and perturbed Sobol' indices via weighted averaging. Since the cost is negligible, it is suggested to compute the optimal perturbation using $\D T_k(\phi)$ and $\D \T_k(\phi)$, and estimate the perturbed Sobol' indices for each perturbation.

Definition~\ref{def:imp} below aids to identify perturbations which change the Sobol' indices but not the relative importance of the variables.

\begin{definition}
\label{def:imp}
Let $\tilde{T}_k$ and $\tilde{\T}_k$ denote the perturbed Sobol' indices and perturbed normalized Sobol' indices, respectively, for some perturbation of $\phi$. The absolute change in the Sobol' indices is
\begin{eqnarray*}
\sum\limits_{k=1}^p \vert T_k - \tilde{T}_k \vert
\end{eqnarray*}
and the relative change in the Sobol' indices is
\begin{eqnarray*}
\sum\limits_{k=1}^p \vert \T_k - \tilde{\T}_k \vert .
\end{eqnarray*}
\end{definition}

It is suggested to consider both the perturbed Sobol' indices which yield the largest absolute and relative changes. This will be further described in Section~\ref{sec:algorithm} and demonstrated in Section~\ref{sec:numerical_results}. We emphasize that the normalized Sobol' indices are a tool to find perturbations changing the relative importance of the variables, but the user will typically use the traditional Sobol' indices to make inferences.

\section{Algorithmic Description}
\label{sec:algorithm}
Algorithm~\ref{alg:robustness} below summarizes our proposed method. In this section, we discuss the user inputs of Algorithm~\ref{alg:robustness} in detail, highlight important algorithmic features of our method, and consider the visualization and interpretation of the results.

It was previously suggested to generate the partition $R_i$, $i=1,2,\dots,M$, with a Regression Tree \cite{breimanbook}. This is a judicious choice because the minimum number of samples in the sets $R_i$ is easily specified. An integer $L$ may be input and the Regression Tree will recursively partition $\Omega$ ensuring that each set of the partition contains at least $L$ samples. This simplicity make Regression Trees attractive in our context. Taking small values of $L$ typically results in $V_M$ being a larger subspace, but will create error when estimating $\D T_u(\phi)\psi_i$ (since there will be fewer samples to estimate the integrals). The determination of $L$ is discussed below. The relationship between $L$ and $M$ depends on the algorithm used to generate the partition; a Regression Tree will uniquely determine $M$ as a function $L$, typically a decreasing function of $L$.

As highlighted in Section~\ref{sec:general}, piecewise constant functions from $V_M$ permit general functional forms for the PDF perturbations. Because the partition size $M$ is constrained by the existing evaluation of $f$, it will typically be coarse by approximation theoretic standards. In order to find the largest possible changes in the Sobol' indices, constrained by the coarseness of the partition, it is advantageous to partition the domain finely in regions where $f$ varies more. This is precisely what a Regression tree, trained to predict $f$, seeks to do. The Regression Tree algorithm recursively partitions the domain, through a greedy algorithm, which minimizes the difference between $f$ and a piecewise constant approximation at each iteration. It begins with the entire domain and refines the partition by considering splits along the directions of the coordinate axes. When $p$ is large it may not split in every coordinate direction. This is acceptable, and in many cases beneficial, as it adapts the partition according to the variability of $f$.

In some cases, as illustrated in Subsection~\ref{sec:lorenz}, partitioning the domain according to the variability of $f$ may result in sets $R_i$, $i=1,2,\dots,M$, where most of the $b_i$'s are small. This is problematic because it limits the size of admissible perturbations. To mitigate this, a Regression Tree may be trained to generate a coarser partition which can be refined by the user to ensure that only a few $b_i$'s are small. We discuss this further in Subsection~\ref{sec:lorenz}.

The norm of the perturbed PDF in \eqref{perturbed_pdf} depends on $\delta$. It was suggested to try various values of $\delta \in [-1,1]$ (equally spaced points in $[-1,1]$) and accept those which meet a convergence tolerance. If Monte Carlo integration is used to estimate the Sobol' indices $\{T_k\}_{k=1}^p$, then the sample standard deviation may be used as a metric for convergence. Let $\sigma_j$ and $\tilde{\sigma}_j$, $j=1,2,\dots,p$, denote the sample standard deviation for the nominal and perturbed Sobol' indices, respectively. For the results presented in this article, the sample standard deviation is estimated by computing the standard derivation of 50 estimates generated by randomly subsampling half of the function evaluations. Assuming that $\sigma_j$, $j=1,2,\dots,p$, are sufficiently small to ensure convergence of the nominal Sobol' indices, it is required that $(\tilde{\sigma}_j/\tilde{T}_j) / (\sigma_j/T_j)$ be less than a threshold. Define
\begin{eqnarray*}
t = \max\limits_{j=1,2,\dots,p} \frac{(\tilde{\sigma}_j/\tilde{T}_j)}{(\sigma_j/T_j)} 
\end{eqnarray*}
and specify a threshold $\tau >1$. The perturbed Sobol' indices are accepted if $t \le \tau$.

The inputs of Algorithm~\ref{alg:robustness} are:
\begin{enumerate}
\item[$\bullet$]  $n$, the number of Monte Carlo samples,
\item[$\bullet$] $L$, the minimum number of samples in each set of the partition, 
\item[$\bullet$] $r$, an integer denoting how many values of $\delta \in [-1,1]$ to consider,
\item[$\bullet$]  and $\tau$, the acceptance threshold for the perturbed Sobol' indices.
\end{enumerate}
The results in Section~\ref{sec:numerical_results} use $L=50$,  $r=60$, and $\tau = 1.5$; the number of Monte Carlo samples required depends on the problem. Numerical evidence, and intuition, indicate that $t$ is approximately a quadratic function of $\delta$ centered at $\delta = 0$. To determine $\delta$, we may solve the scalar nonlinear equation $t(\delta)=\tau$ by evaluating $t(\delta)$ at $r$ equally spaced points in $[-1,1]$. It is not necessary to take large values for $r$; the choice $r=60$ introduces negligible computation and provides sufficient resolution for our purposes. The choice $\tau=1.5$ is considered a reasonable threshold to permit non trivial perturbations without introducing significant numerical errors. Our choice of $L=50$ is the least intuitive of the inputs. To justify this choice, a numerical experiment was performed varying $L=25+5\ell$, $\ell=0,1,\dots,10$. The results, omitted from this article for conciseness, indicate that our method is robust to changes in $L$. If necessary, the user may easily verify the particular choice of inputs used in their application by varying them. The computational cost of this numerical experiment is small.

Lines 2-5 of Algorithm~\ref{alg:robustness} is the Sobol' index estimation and Lines 6-17 is our robustness analysis. In many applications, Line 4 dominates the computational cost and hence the cost of robustness analysis is negligible. Lines 6 and 8 may be done analytically in many applications. The computation in Lines 9-18 is primarily taking sample averages of data on memory so its cost is small. In particular, the nested for loops may appear burdensome, but the operations inside of them are sufficiently simple that they may be executed quickly.

 \begin{algorithm}
\caption{Computation of Sobol' indices with robustness post processing} \label{alg:robustness}
\begin{algorithmic}[1]
\STATE \textbf{Input: } $n$, $L$, $r$, $\tau$
\STATE Draw $n$ samples of $\X$, store them in $X_0 \in \R^{n \times p}$
\STATE Draw $n$ samples of $\X \vert \X_{\sim k}$, store them in $X_k \in \R^{n \times p}$, $k=1,2,\dots,p$
\STATE Evaluate $f(X_j)$, $j=0,1,\dots,p$
\STATE Compute $T_k$, $k=1,2,\dots,p$
\STATE Evaluate $\phi(X_j)$, $j=0,1,\dots,p$
\STATE Generate a partition $\{R_i\}_{i=1}^M$ by using the data $(X_0,f(X_0))$ to train a Regression Tree with a minimum of $L$ data points in each terminal node
\STATE Determine $b_i = \inf_{\x \in R_i} \phi(\x)$, $i=1,2,\dots,M$
\STATE Compute $\D T_k(\phi) \psi_i$, $i=1,2,\dots,M$, $k=1,2,\dots,p$
\STATE Compute $\D \T_k(\phi)\psi_i$, $i=1,2,\dots,M$, $k=1,2,\dots,p$
\FOR{$k$ from 1 to $p$}
\STATE Determine $\psi^{(k,1)} \in V_M$, $\vert \vert \psi^{(k,1)}\vert \vert_V \le 1$, to maximize $\vert \D T_k(\phi) \vert$
\STATE Determine $\psi^{(k,2)} \in V_M$, $\vert \vert \psi^{(k,2)} \vert \vert_V \le 1$, to maximize $\vert \D \T_k(\phi) \vert$
\FOR{$\ell$ from 0 to $r$}
\STATE Compute $\{\tilde{T}_k^{(k,\ell,1)}\}_{k=1}^p$ and $t^{(k,\ell,1)}$ with perturbation $(\phi + \left(-1+\frac{2\ell}{r} \right) \psi^{(k,1)})/N^{(k,\ell,1)}$
\STATE Compute $\{\tilde{T}_k^{(k,\ell,2)}\}_{k=1}^p$ and $t^{(k,\ell,2)}$ with perturbation $(\phi + \left(-1+\frac{2\ell}{r} \right) \psi^{(k,2)})/N^{(k,\ell,2)}$
\ENDFOR
\ENDFOR
\STATE \textbf{Output: } $2p$ sets of perturbed Sobol' indices with largest admissible $t^{(k,\ell,I)} \le \tau$
\STATE \textbf{Note: } $N^{(k,\ell,1)}$, $N^{(k,\ell,2)}$ are constants ensuring the perturbed PDF integrates to one.
\end{algorithmic} 
\end{algorithm}

Algorithm~\ref{alg:robustness} returns a collection of $2p$ sets perturbed Sobol' indices. We suggest extracting the perturbed Sobol' indices with the largest absolute and relative changes to visualize alongside the nominal Sobol' indices, denote them as $\{\tilde{T}_k^a,\tilde{T}_k^r,T_k\}_{k=1}^p$ where the superscripts $a$ and $r$ identify the Sobol' indices with largest absolute and relative changes, respectively. This may be done by querying the collection of perturbed Sobol' indices and creating a bar plot of $\{\tilde{T}_k^a,\tilde{T}_k^r,T_k\}_{k=1}^p$, see Figure~\ref{fig:linear} for an illustration of this. There are several possible scenarios the user may observe:
\begin{enumerate}
\item[$\bullet$] If $\tilde{T}_k^a \approx T_k$, $k=1,2,\dots,p$, then the user may confidently make inferences with the Sobol' indices.
\item[$\bullet$] If $\tilde{T}_k^a \not\approx T_k$, $k=1,2,\dots,p$, but $\tilde{T}_k^r \approx T_k$, $k=1,2,\dots,p$, then the user may confidently make inferences about the relative importance of the variables but not the magnitude of the Sobol' indices.
\item[$\bullet$] If there are variables such that $T_k \approx \tilde{T}_k^a \approx 0$ then they may be considered unimportant.
\item[$\bullet$] If $T_k \approx 0$ but $\tilde{T}_k^a \not\approx 0$ then the user should excise caution treating $x_k$ as unimportant.
\item[$\bullet$] If $T_i > T_j$ but $\tilde{T}_i^r < \tilde{T}_j^r$ then the user may not be certain of the importance of $x_i$ and $x_j$ relative to one another.
\end{enumerate}

 If a particular Sobol' index $T_k$ is of interest, the collection of perturbed Sobol' indices may be queried to asses its robustness. The user may easily visualize all $2p$ of the perturbed indices $\tilde{T}_k$ in a histogram.

\section{Numerical Results}
\label{sec:numerical_results}
In this section, three examples are presented to highlight different properties of our proposed method. The first example analyzes how our robustness analysis changes as more samples are collected. The second example expands on Section~\ref{sec:normalized} by highlighting a case when the largest absolute change in the Sobol' indices yields a small relative change. The final example is an application of our method to the Lorenz system \cite{lorenz63}. We consider two cases in this example to demonstrate the effect of the partitioning on our robustness analysis.

\subsection{g-function Example to Demonstrate Convergence in Samples}
Let
\begin{eqnarray}
\label{gfun}
 f(\X)=\prod\limits_{k=1}^{10} \frac{|4X_k-2|+a_k}{1+a_k},
\end{eqnarray}
where each $X_k$ is independent and uniformly distributed on $[0,1]$, and $a_k =k-1$ for $k=1,2,...,10$. This is the ``g-function" \cite{sobol} commonly used in the GSA literature.

We compute the nominal Sobol' indices and perturbed Sobol' indices of \eqref{gfun}. The number of Monte Carlo samples is varied to analyze the convergence behavior of our robustness estimation, specifically, we use 1,000, 5,000, 10,000, and 50,000 Monte Carlo samples. For each fixed sample size, 32 repetitions of the calculation is performed to understand sampling variability. Figure~\ref{fig:conv} below displays box plots for the estimation of the largest Sobol' index, $T_1$. The center panel is our estimation of $T_1$; the median estimation is nearly constant and the quantiles shrink as the number of samples increases, this reflects convergence of the estimation. The perturbation size $\delta$ is varied between -1 and 1 and it is determined that $\vert \delta \vert = .33$ is the maximum admissible perturbation size for the threshold $\tau=1.5$. The left and right panels show the convergence of $\tilde{T}_1$ with perturbations $\delta = -.33$ and $\delta = .33$, respectively. The shrinking quantiles are very similar to those in the center panel demonstrating that the estimation error in $\tilde{T}_1$ is comparable to the estimation error in $T_1$. The left and right panels have slight decreasing and increasing trends, respectively. This is because the subspace $V_M$ is larger when more samples are taken, thus the perturbations yield larger changes in the Sobol' indices. For this example, the trend is relatively small reflecting the fact that taking a larger subspace does not yield significant changes in the Sobol' index.

\begin{figure}[h]
\centering
\includegraphics[width=.32\textwidth]{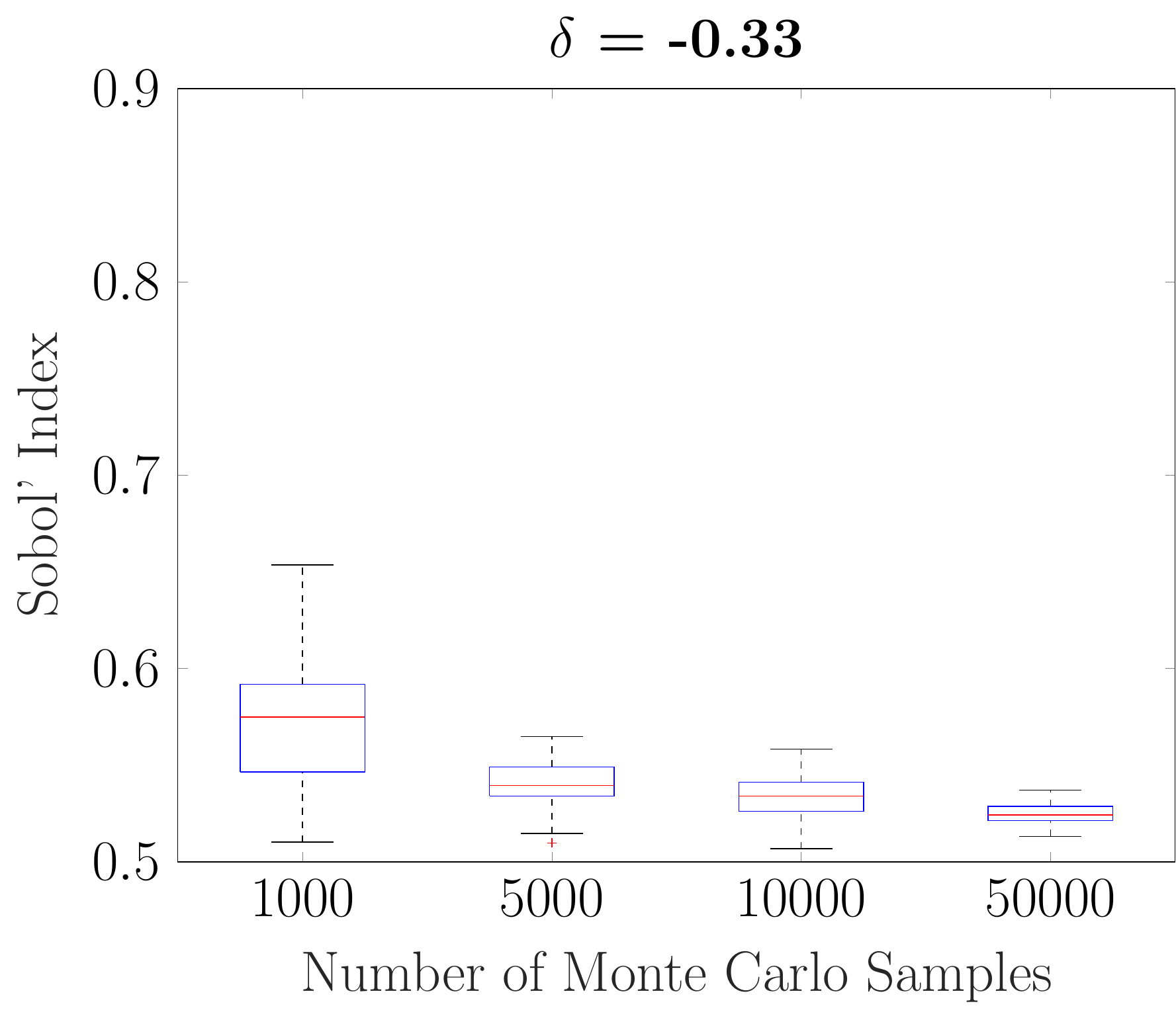}
\includegraphics[width=.32\textwidth]{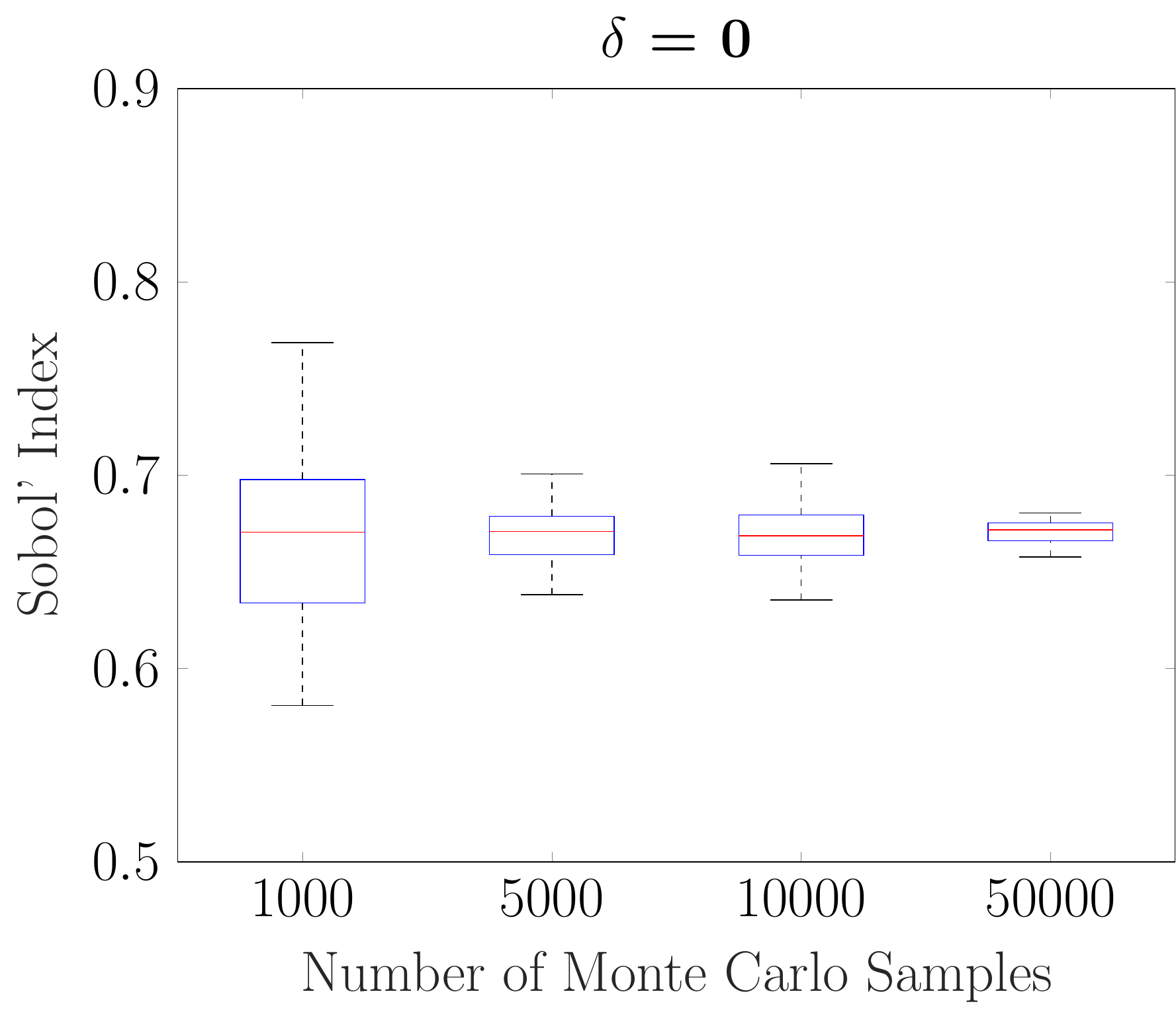}
\includegraphics[width=.32\textwidth]{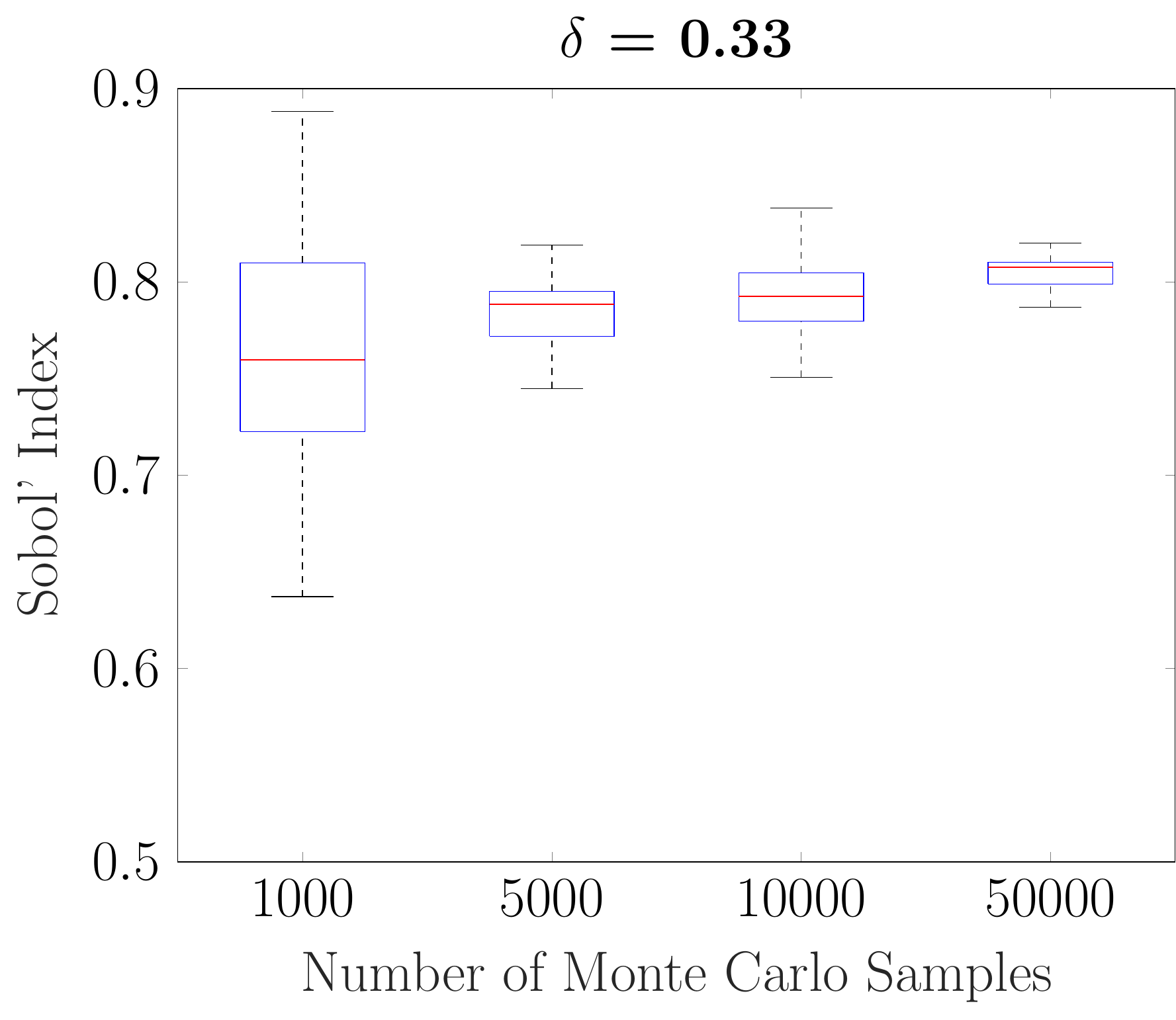}
\caption{Convergence of the Sobol' index of the g-function \eqref{gfun} for variable $x_1$ as the number of Monte Carlo samples vary. Left: perturbed Sobol' index with $\delta=-.33$, center: nominal Sobol' index, right: perturbed Sobol' index with $\delta=.33$. The nominal Sobol' index, computed analytically and rounded, is 0.6743.}
\label{fig:conv}
\end{figure}

\subsection{Linear Example to Demonstrate the Normalized Sobol' Indices}
\label{syn_example}
\begin{figure}[h]
\centering
\includegraphics[width=.49\textwidth]{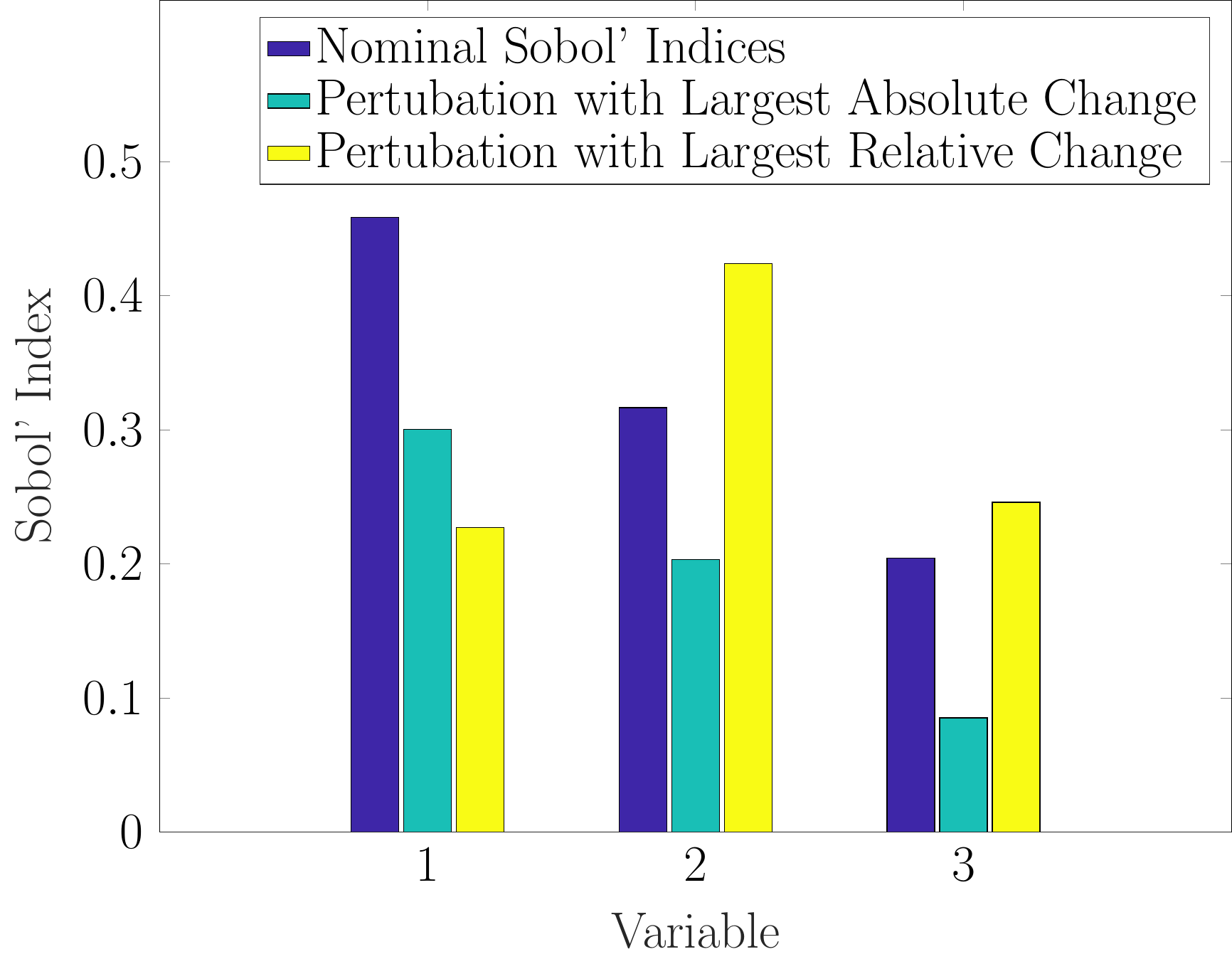}
\caption{Sobol' indices of the linear function \eqref{simpleex} (with independent uniform marginals), the height of each bar indicates the Sobol' index. The blue bars indicate the nominal Sobol' indices; the cyan and yellow bars indicate the Sobol' indices when the PDF of $\x$ was perturbed in extreme cases; cyan: the largest absolute change; yellow: the largest relative change. The nominal Sobol' indices, computed analytically, are $\left( \frac{36}{77},\frac{25}{77},\frac{16}{77}\right)$.}
\label{fig:linear}
\end{figure}

This example illustrates the difference in the largest absolute and relative perturbations of the Sobol' indices. Let $f$ be defined by \eqref{simpleex} and each $X_k$ be independent and uniformly distributed on $[0,1]$, $k=1,2,3$. The Sobol' indices are estimated with 5,000 Monte Carlo samples. Figure~\ref{fig:linear} displays the nominal Sobol' indices of $f$ in blue, the perturbed Sobol' indices with the largest absolute differences change in cyan, and the perturbed Sobol' indices with the largest relative change in yellow.

The largest absolute change of the Sobol' indices corresponds to the case when they are all shifted down but their relative importance does not change. The largest relative change identifies a case where $T_1$ decreases while $T_2$ and $T_3$ increase. The relative importance of the variables change with this perturbation, demonstrating the benefit of considering the largest absolute and relative perturbations.

\subsection{Lorenz System}
\label{sec:lorenz}
This example applies our method to the well known Lorenz system \cite{lorenz63}, a model for atmospheric convection. Sobol' indices were considered for this system in \cite{lorenz_sobol}. The Lorenz system is described by the system of ordinary differential equations
\begin{align*}
& \frac{d y_1}{dt} = \sigma (y_2-y_1) \\
&\frac{d y_2}{dt} = y_1(\rho - y_3) - y_2\\
&\frac{d y_3}{dt} = y_1 y_2 - \beta y_3
\end{align*}
with initial conditions $y_i(0) = \alpha_i$, $i=1,2,3$. Letting $\X = (\sigma, \rho, \beta, \alpha_1,\alpha_2,\alpha_3)$ denote the uncertain parameters, we compute the Sobol' indices of the function
\begin{eqnarray}
\label{lorenz_qoi}
f(\X) = \frac{y_3(1)}{y_2(1)},
\end{eqnarray}
the ratio of the states $y_3$ and $y_2$ at time $t=1$. This choice of $f$ corresponds to a ratio of temperature variations after a duration of 1 time unit.

The distribution of $\X$ is chosen to reflect uncertainty about nominal values of the parameters. Two different cases, in the sub-subsections below, are considered to highlight different features of our method. For each case, $10,000$ Monte Carlo samples are taken for the Sobol' index estimation.

\subsubsection{Lorenz System Case 1}
In this first case we assume the parameters are independent with the uniform distributions given in Table~\ref{tab:case_1_marginals} below. Figure~\ref{fig:lorenz} displays the nominal Sobol' indices in blue, the perturbed Sobol' indices with the largest absolute change in cyan, and the perturbed Sobol' indices with the largest relative change in yellow. Several inferences may be drawn from this result,
\begin{enumerate}
\item[$\bullet$] $\rho$ and $\beta$ are the most influential parameters, although their Sobol' indices and relative importance is uncertain,
\item[$\bullet$] the Sobol' indices for $\sigma$, $\alpha_1$, and $\alpha_2$ and their importance relative to one another is robust,
\item[$\bullet$] $\alpha_3$ has little influence and its small Sobol' index is robust, it may be considered a non-influential parameter.
\end{enumerate}

\begin{table}[h]
\centering
\ra{1.3}
\begin{tabular}{lllll}
\toprule
Parameter & Distribution & Support \\
$\sigma$ & Uniform &  $\left[\frac{97}{10},\frac{103}{10}\right]$ \\
$\rho$ & Uniform  & $\left[\frac{2716}{100},\frac{2884}{100}\right]$ \\
$\beta$ & Uniform  &  $\left[\frac{194}{75},\frac{206}{75}\right]$ \\
$\alpha_1$ & Uniform  & $\left[\frac{4}{5},\frac{6}{5}\right]$ \\
$\alpha_2$ & Uniform  & $\left[\frac{4}{5},\frac{6}{5}\right]$ \\
$\alpha_3$& Uniform  &$\left[\frac{4}{5},\frac{6}{5}\right]$ \\
\bottomrule
\end{tabular}
\caption{Marginal distribution for uncertain parameters in Lorenz system Case 1. The means of $\sigma$, $\rho$, and $\beta$ are the nominal values in \cite{lorenz_sobol}.}
\label{tab:case_1_marginals}
\end{table}

\begin{figure}[h]
\centering
\includegraphics[width=.49\textwidth]{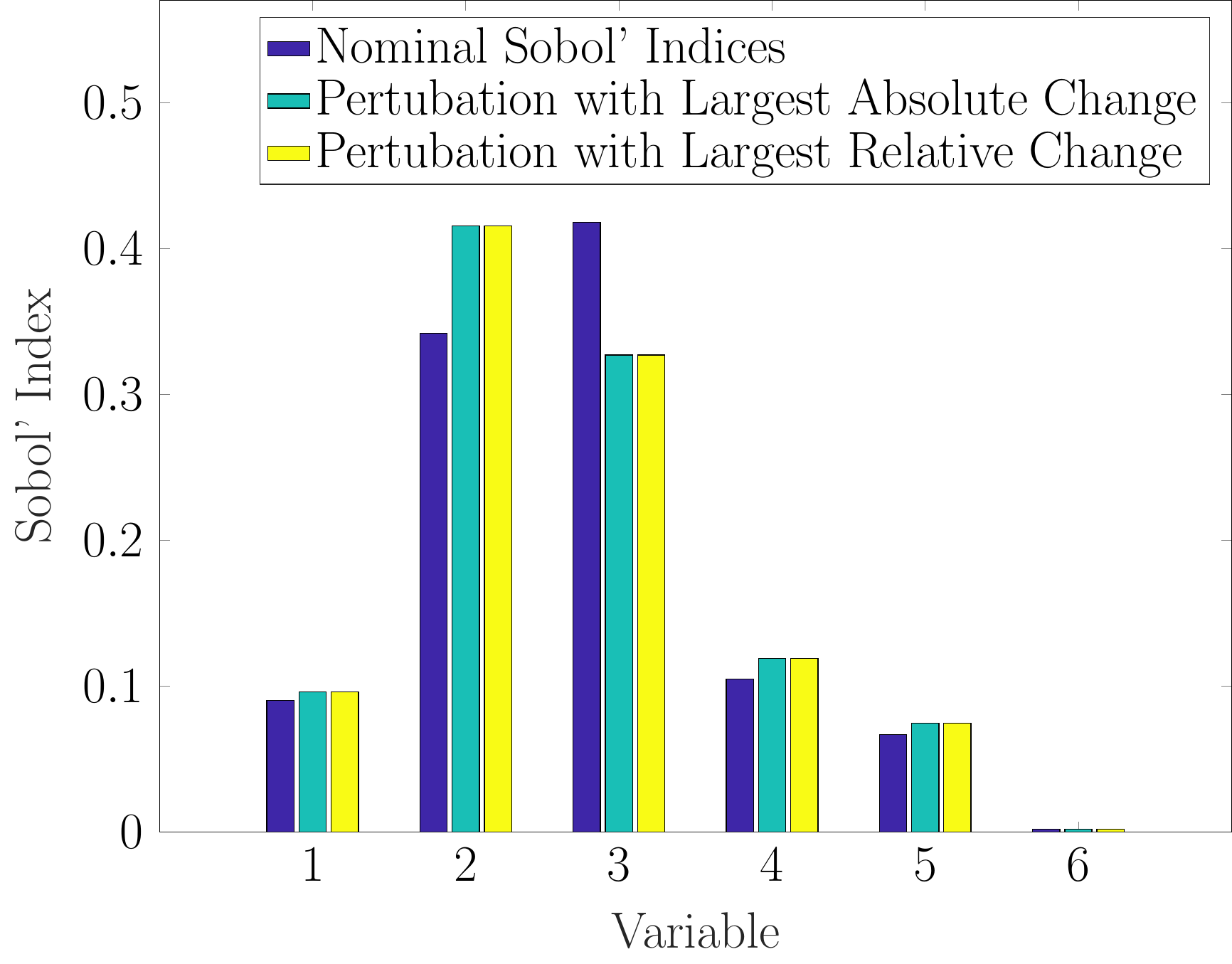}
\caption{Sobol' indices for the Lorenz System \eqref{lorenz_qoi} Case 1 example, the height of each bar indicates the Sobol' index. The blue bars indicate the nominal Sobol' indices; the cyan and yellow bars indicate the Sobol' indices when the PDF of $\x$ was perturbed in extreme cases; cyan: the largest absolute change; yellow: the largest relative change.}
\label{fig:lorenz}
\end{figure}

\subsubsection{Lorenz System Case 2}

In this second case we assume the parameters are independent and that all parameters have the same marginal distribution given in Table~\ref{tab:case_1_marginals} with the exception of $\alpha_3$. Instead of being uniformly distributed on $\left[\frac{4}{5},\frac{6}{5}\right]$ as in Case 1, we take $\alpha_3$ to have a Beta distribution on $\left[\frac{4}{5},\frac{6}{5}\right]$ with shape parameters \footnote{For shape parameters $(a,b)$, a Beta random variable $x$ on $[0,1]$ has PDF $x^a(1-x)^b$.} $(1,4)$. This corresponds to giving greater probability to $\alpha_3 < 1$.

\begin{figure}[h]
\centering
\includegraphics[width=.49\textwidth]{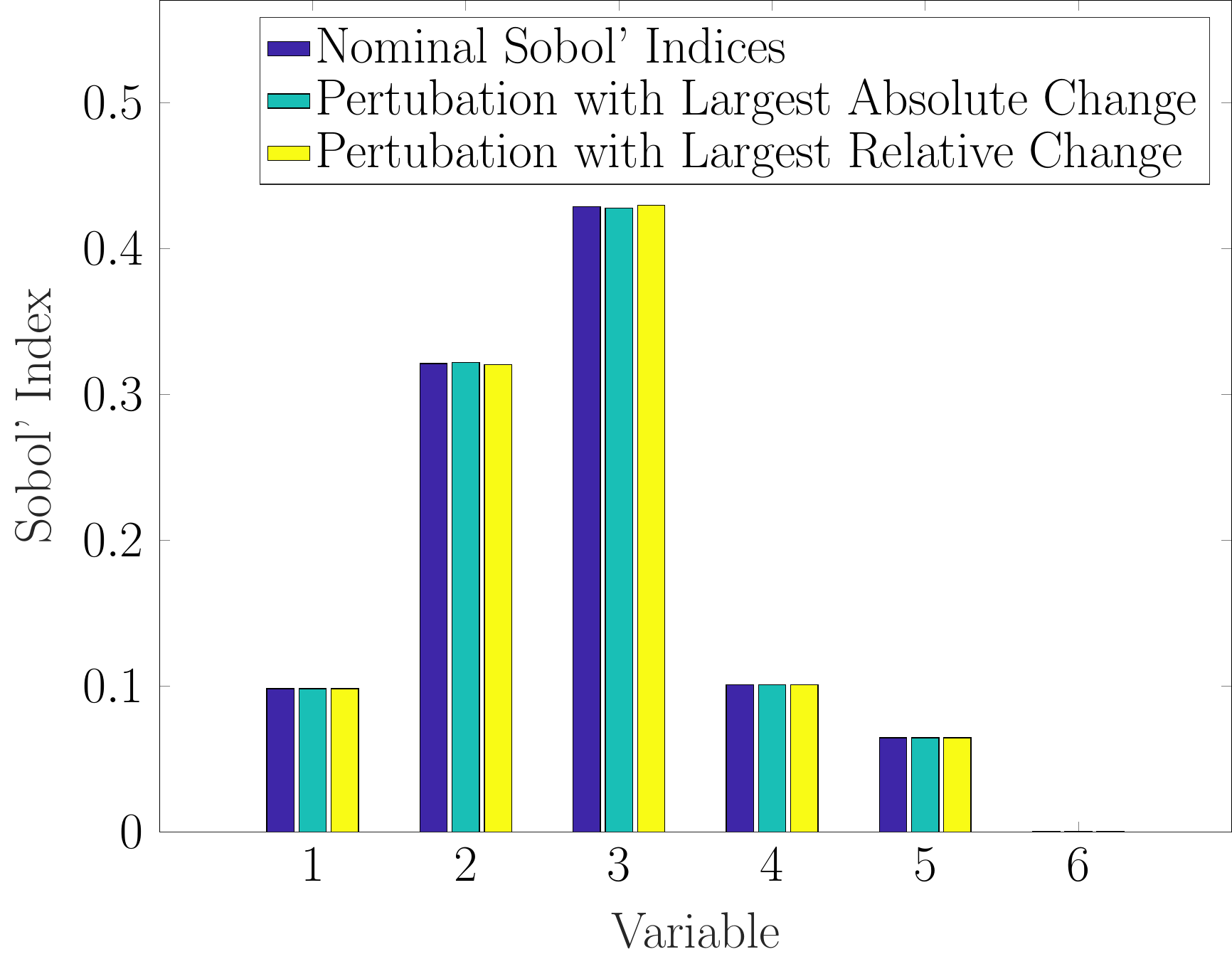}
\includegraphics[width=.49\textwidth]{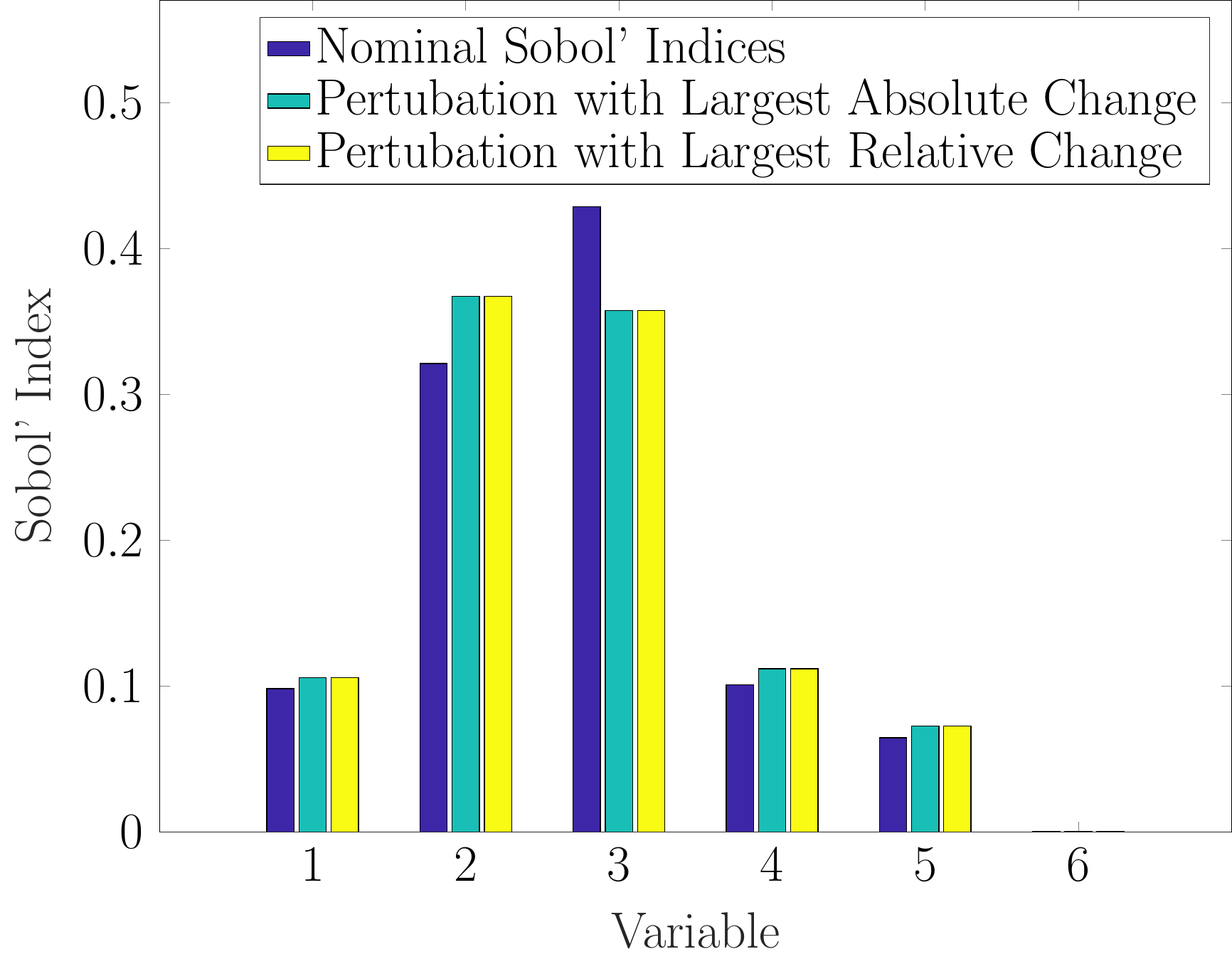}
\caption{Sobol' indices for the Lorenz System \eqref{lorenz_qoi} Case 2 example, the height of each bar indicates the Sobol' index. The blue bars indicate the nominal Sobol' indices; the cyan and yellow bars indicate the Sobol' indices when the PDF of $\x$ was perturbed in extreme cases; cyan: the largest absolute change; yellow: the largest relative change. The left and right panel correspond to generating the partition by training a Regression Tree to: predict $f$ with a minimum of $L$ samples per hyperrectangle (left) and predict $f$ with a minimum of $4L$ samples per hyperrectangle, followed by additional partitioning of $\alpha_3$ (right).}
\label{fig:lorenz_beta}
\end{figure}

A partition is generated by training a Regression Tree to predict $f$. The left panel of Figure~\ref{fig:lorenz_beta} displays the nominal Sobol' indices in blue, the perturbed Sobol' indices with the largest absolute change in cyan, and the perturbed Sobol' indices with the largest relative change in yellow. The results indicate that the Sobol' indices are robust, a different conclusion than was reached in Case 1. This occurs because the Regression Tree never partitioned on $\alpha_3$ so each set $R_i$ contained the entire support of $\alpha_3$. Because the marginal PDF for $\alpha_3$ takes small values on part of its support, namely near $\frac{5}{4}$, the infimum of $\phi$ on each $R_i$ is small. The partition generated by the Regression Tree yielded very small perturbations and as a result did not produce significant changes in the Sobol' indices. 

To alleviate this problem, a partition is generated by a Regression Tree trained to predict $f$ using all of the variables except $\alpha_3$. A minimum of $4L$ samples are requested in each hyperrectangle rather than $L$, as requested previously. This yields a coarser discretization of the other 5 variables. The resulting partition is refined by splitting each set into 4 subsets defined by partitioning at the quantiles of $\alpha_3$. This yields a partition with approximately $L$ samples per subset and a sufficient discretization of $\alpha_3$ to enable larger perturbations. Figure~\ref{fig:lorenz_beta} displays the nominal Sobol' indices in blue, the perturbed Sobol' indices with the largest absolute change in cyan, and the perturbed Sobol' indices with the largest relative change in yellow. Larger changes in the Sobol' indices are observed, as is expected. However the changes are smaller than what was observed in Case 1. This is because the partition used in Case 1 was generated by a Regression Tree which better approximated $f$, and hence allowed for larger perturbations of the Sobol' indices. The general conclusion from this example is that the partition should be generated so that the Regression Tree approximates $f$ as well as possible. If small values of $\phi$ prohibit taking large perturbations, then the partition may be generated with fewer hyperrectangles, followed by a refining of this coarse partition to sufficiently discretize the necessary regions. This may result in a failure to discover the largest possible perturbations, as demonstrated by comparing Case 1 and Case 2.

\section{Conclusion}
\label{sec:conclusion}
This article presents a novel framework in which robustness of the Sobol' indices with respect to the input variables distribution may be assessed. The proposed method permits such analysis to be done at negligible computational cost. For a modeler using Sobol' indices, this robustness analysis can be obtained as a by-product of computing Sobol' indices and may be easily visualized along with the indices themselves. Understanding the robustness of the Sobol' indices to distributional uncertainty prevents the user from making incorrect inferences which have significant consequences. For instance, reducing dimensions by fixing variables with small Sobol' indices---which are not robust---may result in model variations which are not explained in the lower dimensional space.

The method suffers four primary limitations, namely,
\begin{enumerate}
\item the nominal PDF must be compactly supported, \label{l1}
\item perturbations may not change the support of the nominal PDF, \label{l2}
\item the perturbations are only locally optimal, \label{l3}
\item generating a partition is difficult if the distribution of $\X$ is far from being uniform. \label{l4}
\end{enumerate}

The first limitation prohibits a direct application of our method to many commonly used PDF's. This occurs because the Fr\'echet derivative is not well defined if we allow perturbations in the tail of the distribution. However, one may defined a compact subset of the domain where the PDF assumes most of its mass and allow perturbations on this subset while keeping the tail fixed. The compact subset may be chosen large enough that this truncation error is practically irrelevant, for instance, having the probability of the tails less than machine epsilon. This limitation is primarily theoretical and is not a significant practical concern. The greater limitation will be regions of small probability within the compact subset, see the fourth limitation. A theory analogous to what is presented in the article may be developed when the domain is truncated. In practice, the modeler will take a compact subset which contains all of the existing samples and, on the discrete level, the robustness analysis will be identical to what is presented in this article. 

The second limitation occurs because we have formulated the method to work with existing samples. If the support of the PDF increases, then we would need additional evaluations of $f$ in these unexplored regions. 

The third limitation arises because the perturbation direction is determined by maximizing a derivative, which is local. If the PDF to Sobol' index mapping is highly nonlinear, this may not be an adequate. However, finding a globally optimal perturbation requires far more computational effort. A locally optimal perturbation is useful and appropriate, particularly for its computational advantages.


The fourth limitation, as demonstrated in Subsection~\ref{sec:lorenz}, may arise when marginal distributions differ significantly from being uniform. If a small number of marginals do so, this limitation may be mitigated by the approach described in Subsection~\ref{sec:lorenz}. There is ongoing work to develop a variant of our proposed method which removes this limitation by partitioning and taking perturbations on the marginal distributions separately. The method proposed in this article will be most effective when the distribution of $\X$ is approximately uniform, which is a common occurrence in applications where statistical information is poorly known but bounds may be provided from the physics of the problem. 

\section*{Acknowledgement}
The authors thank the two autonomous reviewers for their helpful comments which improved the presentation of this article.

\section*{Appendix}
Proof of Theorem~\ref{thm:frechet}.
\begin{proof}
One may easily observe that $G(\eta)>0$ in a neighborhood of $\phi$ (assuming $f(\x)$ is non constant). It is sufficient to compute the Fr\'echet derivatives of $F$ and $G$, the Fr\'echet derivative of $T_u$ follows from the quotient rule. The Fr\'echet derivatives of 
	\begin{eqnarray*}
	\int_\Omega f(\x) \phi(\x)d\x, \qquad \int_\Omega f(\x)^2 \phi(\x)d\x, \qquad \text{ and } \qquad \int_\Omega \phi(\x)d\x,
	\end{eqnarray*}
	 when considered as operators from $V$ to $\mathbb R$, acting on $\psi$, are easily shown to be
	\begin{eqnarray*}
	\int_\Omega f(\x) \psi(\x)d\x, \qquad \int_\Omega f(\x)^2 \psi(\x)d\x, \qquad \text{ and } \qquad \int_\Omega \psi(\x)d\x,
	\end{eqnarray*}
	respectively, using the definition of the Fr\'echet derivative. The Fr\'echet derivative of $G$ follows from the sum/difference, product, and chain rule of differentiation.
	
	The Fr\'echet derivative of $F$ may be computed by first defining an operator\\ $H:V \to L^\infty(\Omega \times \Omega_u)$,
	\begin{eqnarray*}
	H(\eta) = \eta(\x) \eta(\x') \frac{1}{\int_{\Omega_u} \eta(\x)d\x_u},
	\end{eqnarray*}
	where $\x_{\sim u}'=\x_{\sim u}$. 
	The Fr\'echet derivatives of 
	\begin{eqnarray*}
	\eta(\x), \qquad \eta(\x'), \qquad \text{ and } \qquad \int_{\Omega_u} \eta(\x)d\x_u,
	\end{eqnarray*}
	when considered as operators from $V$ to $L^\infty(\Omega \times \Omega_u)$, acting on $\psi$, are easily shown to be
	\begin{eqnarray*}
	\psi(\x), \qquad \psi(\x'), \qquad \text{ and } \qquad \int_{\Omega_u} \psi(\x)d\x_u,
	\end{eqnarray*}
	respectively, using the definition of the Fr\'echet derivative. The Fr\'echet derivative of $H$ follows from the product and quotient rules of differentiation. The Fr\'echet derivative of $F$ may be easily computed using the Fr\'echet derivative of $H$, the boundedness of $f$, and the chain rule of differentiation.
\end{proof}

\bibliographystyle{siamplain}
\bibliography{PDF_Robustness}

\end{document}

%% file: def.tex
\newcommand{\ra}[1]{\renewcommand{\arraystretch}{#1}}

\newcommand{\R}{\mathbb{R}}

\newcommand{\x}{\mathbf x}

\newcommand{\X}{\mathbf X}

\newcommand{\Va}{\operatorname{Var}}

\newcommand{\T}{\mathscr T}

\newcommand{\E}{\mathbb E}

\newcommand{\D}{\mathscr D}

\newcommand{\LL}{\mathscr L}



%% file: tmp_PDF_Robustness_5_header.tex
\title{Robustness of the Sobol' indices to distributional uncertainty
  \thanks{\funding{This work was supported by the National Science Foundation under grants DMS-1522765 and DMS-1745654.}}}

\author{Joseph Hart
  \thanks{Department of Mathematics, North Carolina State University, Raleigh, NC (\email{jlhart3@ncsu.edu}).}
  \and
  Pierre Gremaud
   \thanks{Department of Mathematics, North Carolina State University, Raleigh, NC (\email{gremaud@ncsu.edu}).}
}

\headers{Robustness of Sobol' indices}
{Joseph Hart and Pierre Gremaud}

%% file: tmp_PDF_Robustness_5_abstract.tex
\begin{abstract}
Global sensitivity analysis (GSA) is used to quantify the influence of uncertain variables in a mathematical model. Prior to performing GSA, the user must specify (or implicitly assume), a probability distribution to model the uncertainty, and possibly statistical dependencies, of the variables. Determining this distribution is challenging in practice as the user has limited and imprecise knowledge of the uncertain variables. This article analyzes the robustness of the Sobol' indices, a commonly used tool in GSA, to changes in the distribution of the uncertain variables. A method for assessing such robustness is developed which requires minimal user specification and no additional evaluations of the model. Theoretical and computational aspects of the method are considered and illustrated through examples.
\end{abstract}

\begin{keywords}
global sensitivity analysis, Sobol' indices, uncertain distributions, deep uncertainty
\end{keywords}

\begin{AMS}
65C60, 62E17
\end{AMS}